\documentstyle[12pt, amsfonts]{amsart}  
\begin{document}     
\def\R{{\mathbb R}}     
\def\Z{{\mathbb Z}}     
\def\C{{\mathbb C}}
\newcommand{\trace}{\rm trace}     
\newcommand{\Ex}{{\mathbb{E}}}     
\newcommand{\Prob}{{\mathbb{P}}}     
\newcommand{\co}{\rm {co\,}}     
\newcommand{\E}{{\mathbb E}}   
\newcommand{\F}{{\cal F}}     
\newtheorem{df}{Definition}     
\newtheorem{theorem}{Theorem}     
\newtheorem{lemma}{Lemma}     
\newtheorem{pr}{Proposition}      
\newtheorem{prob}{Problem}
\newtheorem{cor}{Corollary}   
\newtheorem{remark}{Remark}
\newtheorem{problem}{Problem}
\def\n{\nu}     
\def\sign{\mbox{ sign }}     
\def\a{\alpha}     
\def\N{{\mathbb N}}     
\def\A{{\cal A}}     
\def\L{{\cal L}}
\def\X{{\cal X}}     
\def\F{{\cal F}}     
\def\c{\bar{c}}     
\def\vol{\mbox{\rm Vol}}     
\def\v{\nu}     
\def\d{\delta}     
\def\diam{\mbox{\rm diam}}     
\def\b{\beta}     
\def\t{\theta}     
\def\l{\lambda}     
\def\e{\varepsilon}     
\def\colon{{:}\;}     
\def\pf{\noindent {\bf Proof:  \  }}     
\def\endpf{\begin{flushright} $\Box $\\ \end{flushright}}     
     
\title[Measures of sections of convex bodies]{Measures of sections of convex bodies} 
\author{Alexander Koldobsky}

\address{Department of Mathematics\\ 
University of Missouri\\
Columbia, MO 65211}

\email{koldobskiya@@missouri.edu}

\begin{abstract}  This article is a survey of recent results on slicing inequalities 
for convex bodies. The focus is on the setting of arbitrary measures in place of volume. 

\end{abstract} 
\maketitle

\section{Introduction}

The study of volume of sections of convex bodies is a classical direction in convex geometry.
It is well developed and has numerous applications; see \cite{G3,K4}. 
The question of what happens if volume is replaced by an arbitrary measure on a convex body
has not been considered until very recently, mostly because it is hard to believe
that difficult geometric results can hold in such generality. However, in 2005 Zvavitch \cite{Zv} proved 
that the solution to the Busemann-Petty problem, one of the signature problems in convex geometry,
remains exactly the same if volume is replaced by an arbitrary measure with continuous density. 
It has recently been shown \cite{K6, KM, K8, K9, K10, K11} that several partial results on
the slicing problem, a major open question in the area, can also be extended to arbitrary measures. For example,
it was proved in \cite{K11} that the slicing problem for sections of proportional dimensions has an affirmative answer 
which can be extended to the setting of arbitrary measures. It is not clear yet whether these results 
are representative of something bigger, or it is just an isolated event. We let the reader make the judgement.

\section{The slicing problem for measures}

The slicing problem \cite{Bo1, Bo2, Ba5, MP}, 
asks whether there exists an absolute constant $C$ so that for every origin-symmetric convex 
body $K$ in $\R^n$ of volume 1 there is a hyperplane section of $K$ whose $(n-1)$-dimensional 
volume is greater than $1/C.$
In other words, does there exist a constant $C$ so that for any $n\in \N$ and any
origin-symmetric convex body $K$ in $\R^n$
\begin{equation} \label{hyper}
|K|^{\frac {n-1}n} \le C \max_{\xi \in S^{n-1}} |K\cap \xi^\bot|,
\end{equation}
where  $\xi^\bot$ is the central hyperplane in $\R^n$ perpendicular to $\xi,$ and
$|K|$ stands for volume of proper dimension?
The best current result $C\le O(n^{1/4})$ is due to Klartag \cite{Kla2}, who
slightly improved an earlier estimate of Bourgain \cite{Bo3}.
The answer is known to be affirmative for some special classes of convex bodies,
including unconditional convex bodies (as initially observed by Bourgain; see also \cite{MP, J2,
BN}), unit balls of subspaces of $L_p$  \cite{Ba4, J1, M1}, intersection bodies
\cite[Th.9.4.11]{G3}, zonoids, duals of bodies with bounded volume ratio
\cite{MP}, the Schatten classes \cite{KMP}, $k$-intersection bodies \cite{KPY, K10}.
Other partial results on the problem include \cite{Ba3, BKM, DP, Da, GPV, Kla1,
KlaK, Pa, EK, BaN}; see the book \cite{BGVV} for details.

Iterating (\ref{hyper}) one gets the lower dimensional slicing problem asking whether
the inequality
\begin{equation} \label{lowdimhyper}
|K|^{\frac {n-k}n} \le C^k \max_{H\in Gr_{n-k}} |K\cap H|
\end{equation}
holds with an absolute constant $C,$ where $1\le k \le n-1$ and $Gr_{n-k}$
is the Grassmanian of $(n-k)$-dimensional subspaces of $\R^n.$

Inequality (\ref{lowdimhyper}) was proved in \cite{K11} in the case where $k\ge \lambda n,\ 0<\lambda<1,$
with the constant $C=C(\lambda)$ dependent only on $\lambda.$ Moreover, this was proved in \cite{K11}
for arbitrary measures in place of volume. We consider
the following generalization of the slicing problem to arbitrary measures
and to sections of arbitrary codimension.

\begin{problem} \label{prob}
Does there exist an absolute constant $C$ so that for every $n\in \N,$
every integer $1\le k < n,$ every origin-symmetric convex body $K$ in $\R^n,$ and every measure $\mu$ with non-negative even continuous density $f$ in $\R^n,$
\begin{equation}\label{main-problem}
\mu(K)\ \le\ C^k  \max_{H \in Gr_{n-k}} \mu(K\cap H)\ |K|^{k/n}.
\end{equation}
\end{problem}
Here $\mu(B)=\int_B f$ for every compact set $B$ in $\R^n,$ and $\mu(B\cap H)=\int_{B\cap H} f$
is the result of integration of the restriction of $f$ to $H$ with respect to Lebesgue measure in $H.$
The case of volume corresponds to $f\equiv 1.$

In some cases we will write (\ref{main-problem}) in an equivalent form
\begin{equation}\label{measslicing}
\mu(K)\ \le\ C^k \frac n{n-k}\ c_{n,k} \max_{H \in Gr_{n-k}} \mu(K\cap H)\ |K|^{k/n},
\end{equation}
where $c_{n,k}= |B_2^n|^{\frac {n-k}n}/|B_2^{n-k}|,$ and $B_2^n$ is the unit Euclidean ball in $\R^n.$
It is easy to see that $c_{n,k}\in (e^{-k/2},1)$, and
$ \frac{n}{n-k} \in (1,e^k),$ so these constants
can be incorporated in the constant $C.$

Surprisingly, many partial results on the original slicing problem can be extended 
to the setting of arbitrary measures. Inequality (\ref{main-problem}) holds true in the following cases:
\begin{itemize}
\item for arbitrary $n,K,\mu$ and $k\ge \lambda n,$ where $\lambda\in (0,1),$ with the constant
$C$ dependent only on $\lambda$, \cite{K11};

\item for all $n, K, \mu, k,$ with  $C\le O(\sqrt{n})$, \cite{K8, K9};

\item for intersection bodies $K$ (see definition below), with an absolute constant $C$,
 \cite{K6} for $k=1$, \cite{KM} for all $k;$

\item for the unit balls of $n$-dimensional subspaces of $L_p,\ p>2,$
with $C\le O(n^{1/2-1/p}),$ \cite{K10};

\item for the unit balls of $n$-dimensional normed spaces that embed in $L_p,\ p\in (-n, 2]$, 
with $C$ depending only on $p,$ \cite{K10}; 

\item for unconditional convex bodies, with an absolute constant $C,$ \cite{K11};

\item for duals of convex bodies with bounded volume ratio, with an absolute constant $C,$ \cite{K11};

\item for $k=1$ and log-concave measures $\mu$, with $C\le O(n^{1/4}),$ \cite{KZ}.
\end{itemize}

The proofs of these results are based on stability in volume
comparison problems introduced in \cite{K5} and developed in \cite{K6, KM, K7, K8,
K9, K10,K13}. Stability reduces Problem \ref{prob} to estimating the outer 
volume ratio distance from a convex body to the classes of generalized intersection bodies.
The concept of an intersection body was introduced by Lutwak \cite{Lu} in connection with the
Busemann-Petty problem.

A closed bounded set $K$ in $\R^n$ is called a {star body}  if 
every straight line passing through the origin crosses the boundary of $K$ 
at exactly two points different from the origin, the origin is an interior point of $K,$
and the boundary of $K$ is continuous.

For $1\le k \le n-1,$ the classes ${\cal{BP}}_k^n$ of generalized $k$-intersection bodies in $\R^n$ 
were introduced by Zhang \cite{Z3}. The case $k=1$ represents the original class of intersection bodies ${\cal{I}}_n={\cal{BP}}_1^n$ of Lutwak \cite{Lu}. We define ${\cal{BP}}_k^n$ as the closure in the radial metric
of radial $k$-sums of finite collections of origin-symmetric ellipsoids (the equivalence of this definition to the
original definitions of Lutwak and Zhang was established by Goodey and Weil \cite{GW} for $k=1$
and by Grinberg and Zhang \cite{GrZ} for arbitrary $k.)$ Recall that the radial $k$-sum of star bodies
$K$ and $L$ in $\R^n$ is a new star body $K+_kL$ whose radius in every direction $\xi\in S^{n-1}$ 
is given by
$$r^k_{K+_kL}(\xi)= r^k_{K}(\xi) + r^k_{L}(\xi).$$
The radial metric in the class of origin-symmetric star bodies is defined by
$$\rho(K,L)=\sup_{\xi\in S^{n-1}} |r_K(\xi)-r_L(\xi)|.$$
The following stability theorem was proved in \cite{K11} (see \cite{K6, KM} for slightly 
different versions).

\begin{theorem}\label{stab2}(\cite{K11})
Suppose that $1\le k \le n-1,$ $K$ is a generalized $k$-intersection body in $\R^n,$  $f$
is an even continuous non-negative function on $K,$ and $\e>0.$ If
$$
\int_{K\cap H} f \ \le \e,\qquad \forall H\in Gr_{n-k},
$$
then
$$
\int_K f\ \le \frac n{n-k}\ c_{n,k}\ |K|^{k/n}\e.
$$
The constant is the best possible. Recall that $c_{n,k}\in (e^{-k/2},1).$
\end{theorem}
Define the outer volume ratio distance from an origin-symmetric star body $K$ in $\R^n$ to the class 
${\cal{BP}}_k^n$ of generalized $k$-intersection bodies by
$${\rm {o.v.r.}}(K,{\cal{BP}}_k^n) = \inf \left\{ \left( \frac {|D|}{|K|}\right)^{1/n}:\ K\subset D,\ D\in {\cal{BP}}_k^n \right\}.$$

Theorem \ref{stab2} immediately implies a slicing inequality for arbitrary measures and origin-symmetric star bodies.

\begin{cor} \label{lowdim} Let $K$ be an origin-symmetric star body in $\R^n.$ Then for any measure $\mu$
with even continuous density on $K$ we have
$$\mu(K)\le  \left({\rm{ o.v.r.}}(K,{{\cal{BP}}_k^n})\right)^k \frac n{n-k}\ c_{n,k} \max_{H\in Gr_{n-k}} \mu(K\cap H)\ |K|^{k/n}.$$
\end{cor}

Thus, stability reduces Problem 1 to estimating the outer volume ratio distance from $K$ to the class
of generalized $k$-intersection bodies. The results on Problem \ref{prob} mentioned above were all obtained by
estimating this distance by means of various techniques from the local theory of Banach spaces.
For example, the solution to the slicing problem for sections of proportional dimensions follows from
an estimate obtained in \cite{KPZ}: for any origin-symmetric convex 
body $K$ in $\R^n$ and any $1\le k \le n-1,$ 
\begin{equation}\label{kpz}
{\rm o.v.r.}(K,{\cal{BP}}_k^n) \le C_0 \sqrt{\frac{n}{k}}\left(\log\left(\frac{en}{k}\right)\right)^{3/2},
\end{equation}
where $C_0$ is an absolute constant.
The proof of this estimate in \cite{KPZ} is quite involved. It uses covering numbers, Pisier's generalization
of Milman's reverse Brunn-Minkowski inequality, properties of intersection bodies.
Combining this with Corollary \ref{lowdim}, one gets
\begin{theorem}\label{proport} (\cite{K11}) If the codimension of sections $k$ satisfies $\lambda n\le k$ for some $\lambda\in (0,1),$ 
then for every origin-symmetric convex body $K$ in $\R^n$ and every measure $\mu$ with continuous
non-negative density in $\R^n,$
$$ \mu(K)\ \le\  C^k 
\left(\sqrt{\frac{(1-\log \lambda)^3}{\lambda}}\right)^k 
 \max_{H \in Gr_{n-k}} \mu(K\cap H)\ |K|^{k/n},$$
where $C$ is an absolute constant.
\end{theorem}
For arbitrary $K,\mu$ and $k$ the best result so far is the following $\sqrt{n}$ estimate; see \cite{K8, K9}.
By John's theorem, for any origin-symmetric convex body $K$ there exists an ellipsoid
${\cal{E}}$ so that $ \frac 1{\sqrt{n}} {\cal{E}}\subset K \subset {\cal{E}}.$ Since every ellipsoid
is a generalized $k$-intersection body for every $k,$ we get that
$${\rm o.v.r.}(K,{\cal{BP}}_k^n) \le \sqrt{n}.$$
By Corollary \ref{lowdim},
$$
\mu(K)\ \le\  n^{k/2} \frac n{n-k}\ c_{n,k}\max_{H\in Gr_{n-k}}  \mu(K\cap H)\ |K|^{k/n}.
$$

Note that the condition that the measure $\mu$ has continuous density is necessary in Problem \ref{prob}.
A discrete version of inequality (\ref{main-problem}) was very recently established (with a constant depending only
on the dimension) in \cite{AHZ}.

\section{The isomorphic Busemann-Petty problem} In 1956, Busemann and Petty \cite{BP}
asked the following question. Let $K,L$ be origin-symmetric convex bodies 
in $\R^n$ such that
\begin{equation}\label{bp-condition}
\left|K\cap\xi^\perp\right| \le \left|L\cap\xi^\perp\right|,\qquad \forall \xi\in S^{n-1}.
\end{equation}
Does it necessarily follow that $\left|K\right| \le \left|L\right| ?$
The problem was solved at the end of the 1990's in 
a sequence of papers \cite{LR, Ba1, Gi, Bo4, 
Lu, P, G1, G2, Z1, K1, K2, Z2,
GKS}; see \cite[p.3]{K4} or \cite[p.343]{G3} for the solution and its history.
The answer is affirmative if $n\le 4$, and it is negative if $n\ge 5.$

The {lower dimensional Busemann-Petty problem} asks the same question for
sections of lower dimensions. Suppose that $1\le k \le n-1,$ and $K,L$ 
are origin-symmetric convex bodies in $\R^n$ such that
\begin{equation}\label{ldbp-condition}
|K\cap H|\le |L\cap H|,\qquad \forall H\in Gr_{n-k}.
\end{equation}
Does it follow that $|K|\le |L|?$ It was proved in \cite{BZ} (see also \cite{K3, K4, RZ,
M2} for different proofs) that the answer is negative if the dimension of sections $n-k>3.$
The problem is still open for two- and three-dimensional sections ($n-k=2,3,\ n\ge 5).$

Since the answer to the Busemann-Petty problem is negative in most dimensions, it makes sense to ask 
whether the inequality for volumes holds up to an absolute constant, namely, does there exist an absolute 
constant $C$ such that inequalities (\ref{bp-condition}) imply $|K|\le C\ |L|\ ?$ 
This question is known as
the isomorphic Busemann-Petty problem, and in the hyperplane case it is equivalent
to the slicing problem; see \cite{MP}.
A version of this problem for sections of proportional 
dimensions was proved in \cite{K12}. 
\begin{theorem}\label{bp-proport}(\cite{K12})
Suppose that $0<\lambda<1,$ $k>\lambda n,$ and $K,L$ are
origin-symmetric convex bodies in $\R^n$ satisfying the inequalities
$$|K\cap H|\le |L\cap H|,\qquad \forall H\in Gr_{n-k}.$$
Then
$$|K|^{\frac{n-k}n}\le \left(C(\lambda)\right)^k|L|^{\frac{n-k}n},$$
where $C(\lambda)$ depends only on $\lambda.$ 
\end{theorem}
This result implies Theorem \ref{proport} in the case of volume. It is not clear, however,
whether Theorem \ref{proport} can be directly used to prove Theorem \ref{bp-proport}.
\smallbreak
Zvavitch \cite{Zv} has found a remarkable generalization of the Busemann-Petty problem
to arbitrary measures in place of volume. Suppose that $1\le k <n,$ $\mu$ is a measure with 
even continuous density $f$ in $\R^n,$ and $K$ and $L$ are origin-symmetric convex bodies
in $\R^n$ so that
\begin{equation}\label{bp-measure}
\mu(K\cap \xi^\bot) \le \mu(L\cap \xi^\bot), \qquad  \forall \xi\in S^{n-1}.
\end{equation}
Does it necessarily follow that $\mu(K)\le \mu(L)?$ The answer is the same as for volume -
affirmative if $n\le 4$ and negative if $n\ge 5.$  An isomorphic version was recently proved in \cite{KZ}, namely,
for every dimension $n$ inequalities (\ref{bp-measure}) imply $\mu(K)\le \sqrt{n}\ \mu(L).$ It is not
known whether the constant $\sqrt{n}$ is optimal for arbitrary measures. Also there is no known direct connection between the isomorphic Busemann-Petty problem for
arbitrary measures and Problem \ref{prob}.

\section{Projections of convex bodies.} The projection analog of the Busemann-Petty problem 
is known as Shephard's problem, posed in 1964 in \cite{Sh}.
Denote by $K\vert \xi^\bot$ the orthogonal projection of $K$ to $\xi^\bot.$
Suppose that $K$ and $L$ are origin-symmetric convex bodies in 
$\R^n$ so that  $|K\vert \xi^\bot|\le |L\vert \xi^\bot|$  for every  $\xi\in S^{n-1}.$ Does it follow that
$|K|\le |L|?$ The problem was solved by Petty \cite{Pe} and 
Schneider \cite{Sch}, independently, and 
the answer if affirmative only in dimension 2. 

Both solutions use the fact that
the answer to Shephard's problem is affirmative in every dimension 
under the additional assumption that $L$ is a projection body.
An origin symmetric convex body $L$ in $\R^n$ 
is called a {projection body} if there exists another convex body $K$ 
so that the support function of $L$ in every direction is equal to
the volume of the hyperplane projection of $K$ to this direction: 
for every $\xi\in S^{n-1},$ 
$$
h_{L}(\xi) = |K\vert\xi^{\bot}|.
$$
The support function $h_L(\xi)=\max_{x\in L} |(\xi, x)|$ is equal to the
dual norm $\|\xi\|_{L^*},$ where $L^*$ denotes the polar body of $L.$

Separation in Shephard's problem was proved in \cite{K5}. 
\begin{theorem} \label{main-proj1} (\cite{K5}) 
Suppose that $\e>0$,  $K$ and $L$ are origin-symmetric
convex bodies in $\R^n,$ and $L$ is a projection body.  If 
$|K\vert \xi^\bot|\le |L\vert \xi^\bot| - \e$ for every $\xi\in S^{n-1},$
then
$|K|^{\frac{n-1}n}  \le |L|^{\frac{n-1}n} - c_{n,1} \e,$
where $c_{n,1}$ is the same constant as in Theorem \ref{stab2}; recall that $c_{n,1}>1/\sqrt{e}.$
\end{theorem}

Stability in Shephard's problem turned out to be more difficult, and it was proved in \cite{K13} only up to 
a logarithmic term
and under an additional assumtpion that the body $L$ is isotropic. Recall that a convex body $D$ in 
$\R^n$ is isotropic if $|D|=1$ and $\int_{D} (x,\xi)^2 dx$ is a constant function
of $\xi\in S^{n-1}.$ Every convex body has a linear image that is isotropic; see [BGVV].
\begin{theorem} \label{stabproj} (\cite{K13}) 
Suppose that $\e>0$,  $K$ and $L$ are origin-symmetric
convex bodies in $\R^n,$ and $L$ is a projection body which is a dilate of an isotropic body.  
If $|K\vert \xi^\bot|\le |L\vert \xi^\bot| + \e$ for every $\xi\in S^{n-1},$
then
$|K|^{\frac{n-1}n}  \le |L|^{\frac{n-1}n} + C \e \log^2n,$
where $C$ is an absolute constant. 
\end{theorem}
The proof is based on an estimate for the mean width 
of a convex body obtained by E.Milman [M3].

The projection analog of the slicing problem reads as
$$
|K|^{\frac {n-1}n} \ge c \min_{\xi\in S^{n-1}} |K\vert \xi^\bot|,
$$
and it was solved by Ball [Ba2], who proved that 
$c$ may and has to be of the order $1/\sqrt{n}.$

The possibility of extension of Shephard's problem and related stability and separation
results to arbitrary measures is an open question.
Also, the lower dimensional Shephard problem was solved by Goodey and Zhang [GZ],
but stability and separation for the lower dimensional case have not been established.
 
Stability and separation for projections have an interesting application to 
surface area. If $L$ is a projection body, so is $L+\e B_2^n$ for every $\e>0.$ 
Applying stability in Shephard's problem to this 
pair of bodies, dividing by $\e$ and sending $\e$ to zero, one gets a hyperplane 
inequality for surface area (see \cite{K7}): if $L$ is a projection body, then 
\begin{equation}\label{surf-proj-min}
S(L)\ge c \min_{\xi\in S^{n-1}} S(L\vert \xi^\bot)\ |L|^{\frac 1n}.
\end{equation}
On the other hand, applying separation to any projection body $L$ which is a dilate of a body in
isotropic position (see \cite{K13})
\begin{equation}\label{surf-proj-max}
S(L)\le C\log^2n \max_{\xi\in S^{n-1}} S(L\vert \xi^\bot)\ |L|^{\frac 1n}.
\end{equation}
Here $c$ and $C$ are absolute constants, and $S(L)$ is surface area.

\end{document}